\documentclass[12pt]{amsart}
\usepackage{amssymb}
\usepackage{amsmath}
\usepackage{amsthm}
\newcommand{\E}{{\mathcal E}}
\renewcommand{\H}{{\mathcal H}}
\newcommand{\I}{{\mathcal I}}
\newcommand{\N}{{\mathcal N}}
\newcommand{\C}{\ensuremath{\mathbb{C}}}
\newcommand{\R}{\ensuremath{\mathbb{R}}}
\newcommand{\vk}{\varkappa}
\newcommand{\p}{\partial}
\newcommand{\s}{{\rm Symb}}
\renewcommand{\S}{\Sigma}
\newcommand{\Vect}{{\rm Vect}}
\newcommand{\T}{{T^\ast M}}
\newcommand{\D}{{\mathcal D}}
\newcommand{\W}{{\mathcal W}}
\newtheorem{lemma}{Lemma}

\newtheorem{proposition}{Proposition}
\newtheorem{theorem}{Theorem}
\newtheorem{corollary}{Corollary}
\begin{document}

\title[Fedosov's Formal Symplectic Groupoids]{Fedosov's Formal Symplectic Groupoids and Contravariant Connections}

\author[A.V. Karabegov]{Alexander V. Karabegov}
\address[Alexander V. Karabegov]{Department of Mathematics and Computer Science, Abilene Christian University, ACU Box 28012, 252 Foster Science Building, Abilene, TX 79699-8012}
\email{axk02d@acu.edu}

\begin{abstract}
  Using Fedosov's approach we give a geometric construction of a formal symplectic groupoid over any Poisson manifold endowed with a torsion-free Poisson contravariant connection. In the case of K\"ahler-Poisson manifolds this construction provides, in particular, the formal symplectic groupoids with separation of variables. We show that the dual of a semisimple Lie algebra does not admit torsion-free Poisson contravariant connections.
\end{abstract}

\maketitle

\medskip

\section{Introduction}

\bigskip

A symplectic groupoid over a Poisson manifold $M$ is a symplectic manifold $\S$ endowed with a partially defined multiplication and the source, target, inverse, and unit mappings satisfying several axioms. In particular, the source and the target mappings are a Poisson and an anti-Poisson mappings from $\S$ to $M$, respectively. Symplectic groupoids play the r{\^o}le of semiclassical counterparts of associative algebras treated as quantum objects. Symplectic groupoids were introduced independently by Karasev  \cite{Ka}, Weinstein \cite{W}, \cite{CDW}, and Zakrzewski \cite{Z}. There is a corresponding notion of a formal symplectic groupoid on the formal neighborhood $(\S,\Lambda)$ of a Lagrangian submanifold $\Lambda$ of a symplectic manifold $\S$ whose principal example is the formal neighborhood $(\S,\Lambda)$ of the Lagrangian unit space $\Lambda$ of a symplectic groupoid on a symplectic manifold $\S$ (see \cite{FSG}). Formal symplectic groupoids were first introduced in \cite{CDF} in terms of formal generating functions of their (formal) Lagrangian product spaces. It was shown in \cite{FSG} that to each natural deformation quantization on a Poisson manifold $M$  there corresponds a canonical formal symplectic groupoid on $(\T,Z)$, where $Z$ is the zero section of $\T$. The main result of \cite{CDF} is the description of the formal symplectic groupoid of Kontsevich deformation quantization.  The formal symplectic groupoid of Fedosov's star-product was described in \cite{Deq}. This paper is motivated by the following observation. On the one hand, it is known that deformation quantizations with separation of variables (also known as deformation quantizations of the Wick type, see \cite{CMP1} and \cite{BW}) are a particular case of Fedosov's deformation quantizations (see \cite{N}). On the other hand, it was shown in \cite{FSG} that the corresponding formal symplectic groupoids `with separation of variables' can be naturally extended from K\"ahler manifolds to K\"ahler-Poisson manifolds, while it is impossible to extend the star-products with separation of variables to the K\"ahler-Poisson manifolds in a naive direct way (see \cite{Third}). In this paper we show that the construction  of the formal symplectic groupoids of Fedosov's deformation quantizations from \cite{Deq} can be naturally extended to the Poisson manifolds endowed with a torsion-free Poisson contravariant connection. We call the formal symplectic groupoids obtained via this construction Fedosov's formal symplectic groupoids. 

On a K\"ahler-Poisson manifold, there is a natural torsion-free Poisson contravariant connection which we call the K\"ahler-Poisson contravariant connection. We show that Fedosov's formal symplectic groupoid constructed with the use of the K\"ahler-Poisson contravariant connection is a formal symplectic groupoid with separation of variables.

Any symplectic manifold admits symplectic (torsion-free) connections and therefore Poisson torsion-free contravariant connections. However, this is not the case for general Poisson manifolds. We prove that the dual space of a semisimple Lie algebra does not admit a torsion-free Poisson contravariant connection.

{\bf Acknowledgments.} I am very grateful to Professor Alan Weinstein for giving me the opportunity to present this work at Northern California Symplectic Geometry Seminar (NCSGS) at Stanford University. I would like to thank Victor Ginzburg and other participants of the NCSGS for interesting discussions.

\medskip

\section{Linear Contravariant Connections}\label{S:linear}

\bigskip

Contravariant derivatives were introduced by I. Vaisman in \cite{Vais}.  The corresponding notion of a contravariant connection was extensively studied by R. L. Fernandes in \cite{RLF}.

Let $M$ be a Poisson manifold endowed with the Poisson bivector field $\Pi$. Then the Poisson bracket of functions $f,g \in C^\infty(M)$ is given by
\[
             \{f,g\} = \Pi(df,dg).
\]
Define a bundle map $\#: T^*M \to TM$ by the formula
\[
   \langle \beta, \#\alpha\rangle = \Pi(\alpha, \beta),
\]
where $\alpha, \beta \in \Omega^1(M)$ are 1-forms on $M$ and $\langle \cdot,\cdot \rangle$ is the natural pairing of $T^*M$ and $TM$. It is known that on the space $\Omega^1(M)$ of 1-forms on $M$ there is a Lie bracket
\[
    [\alpha,\beta] = {\mathcal L}_{\#\alpha} \beta - {\mathcal L}_{\#\beta} \alpha - d\Pi(\alpha,\beta),
\]
where $\alpha, \beta \in \Omega^1(M)$ and ${\mathcal L}$ denotes the Lie derivative.
If $\alpha = df$ and $\beta = dg$ for $f,g \in C^\infty(M)$, then
\begin{equation}\label{E:exact}
       [df,dg] = d\{f,g\}.
\end{equation}
The mapping $\#$ induces a Lie algebra homomorphism from $\Omega^1(M)$ to the Lie algebra of vector fields on $M$, ${\rm Vect}(M)$, so that for $\alpha, \beta \in \Omega^1(M)$
\begin{equation}\label{E:hom}
         \#[\alpha,\beta] = [\#\alpha,\#\beta].
\end{equation}

A contravariant connection $\nabla^\bullet$ on a vector bundle $E$ on $M$ is a bilinear mapping $\nabla^\bullet: \Gamma(T^*M \otimes E) \to \Gamma(E)$ satisfying the axioms 

\medskip

 i) $\nabla^{f \alpha} s = f \nabla^\alpha s$;

 ii) $\nabla^\alpha (fs) = f \nabla^\alpha s + \#\alpha(f) s$,\\
where $f\in C^\infty(M),\ \alpha \in \Omega^1(M),$ and $s \in \Gamma(E)$.

A covariant connection $\nabla_\bullet$ on $E$ induces a contravariant connection $\nabla^\bullet$ as follows:
\begin{equation}\label{E:ind}
    \nabla^\alpha = \nabla_{\#\alpha}.
\end{equation}
It is known that there exist contravariant connections that are not induced by a covariant connection.

Let $\nabla^\bullet$ be a contravariant connection on $T^*M$. Its torsion $T^\nabla$ and curvature $R^\nabla$ are defined by the formulas
\[
   T^\nabla(\alpha, \beta) = \nabla^\alpha \beta - \nabla^\beta \alpha - [\alpha, \beta]
\]
and 
\[
   R^\nabla(\alpha,\beta)\gamma = \nabla^\alpha \nabla^\beta \gamma - \nabla^\beta \nabla^\alpha \gamma - \nabla^{[\alpha,\beta]}\gamma,
\]
where $\alpha,\beta,\gamma$ are 1-forms on $M$. The transposed contravariant connection $^t\nabla^\bullet$ of a contravariant connection $\nabla^\bullet$ is defined as follows:
\[
        ^t\nabla^\alpha \beta = \nabla^\beta \alpha + [\alpha, \beta],
\]
so that $T^\nabla(\alpha,\beta) = \nabla^\alpha \beta - ^t\!\nabla^\alpha \beta$. The transposed connection of $^t\nabla^\bullet$ is $\nabla^\bullet$.

To any linear covariant connection $\nabla_\bullet$ on $TM$ there corresponds a covariant connection on $T^*M$ which will be denoted also by $\nabla_\bullet$. In local coordinates $\{x^i\}$ 
\[
\nabla_{\p_i} \p_j = \Gamma^k_{ij} \p_k \mbox{ and } \nabla_{\p_i} dx^j = - \Gamma^j_{ik} dx^k,
\]
where $\p_i = \p/\p x^i$ and $\Gamma^k_{ij}$ are the Christoffel symbols of $\nabla_\bullet$.
Similarly, to any contravariant connection $\nabla^\bullet$ on $TM$ there corresponds a contravariant connection on $T^*M$ which will be denoted also by $\nabla^\bullet$. In local coordinates 
\[
\nabla^{dx^i} \p_j =  \Gamma_j ^{ik} \p_k \mbox{ and } \nabla^{dx^i} dx^j = - \Gamma_k^{ij} dx^k,
\]
where $\Gamma_k^{ij}$ are the Christoffel symbols of $\nabla^\bullet$. In particular, if a contravariant connection $\nabla^\bullet$ is induced by a covariant connection $\nabla_\bullet$ according to (\ref{E:ind}), then
\[
          \Gamma^{ij}_k = \pi^{il} \Gamma^j_{lk},
\]
where $\pi^{ij}$ is the Poisson tensor corresponding to $\Pi$.

To a given a covariant connection $\nabla_\bullet$ on $TM$ there corresponds the transposed connection $^t \nabla_\bullet$ on $TM$ such that $^t \Gamma^k_{ij} = \Gamma^k_{ji}$, where $\Gamma^k_{ij}$ and $^t \Gamma^k_{ij}$ are the Christoffel symbols of $\nabla_\bullet$ and $^t \nabla_\bullet$, respectively. For vector fields $X,Y \in {\rm Vect}(M)$
\begin{equation}\label{E:transp}
   ^t \nabla_X Y = \nabla_Y X + [X,Y],  
\end{equation}
so that the torsion $T_\nabla(X,Y)$ of $\nabla_\bullet$ is given by the formula
\[
    T_\nabla(X,Y) = \nabla_X Y - ^t\!\nabla_X Y.
\]
The transposed connection of $^t \nabla_\bullet$ is $\nabla_\bullet$. Given a covariant connection $\nabla_\bullet$, one can construct a pair of contravariant connections $(\nabla^\bullet, ^\dagger\!\nabla^\bullet)$ induced by $\nabla_\bullet$ and $^t \nabla_\bullet$, respectively, according to (\ref{E:ind}). The connection $^\dagger \nabla^\bullet$ is in general different from the transposed contravariant connection $^t \nabla^\bullet$. For 1-forms $\alpha,\beta \in \Omega^1(M)$ it follows from (\ref{E:hom}), (\ref{E:ind}), and (\ref{E:transp}) that
\begin{equation}\label{E:agree}
      ^\dagger \nabla^\alpha \# \beta - \nabla^\beta \# \alpha - \# [\alpha,\beta] =0.
\end{equation}
A pair $(\nabla^\bullet, ^\dagger\! \nabla^\bullet)$ of contravariant connections satisfying (\ref{E:agree}) is a natural counterpart of a single covariant connection. We will call such connections {\it associated}. In general, one cannot recover one of the associated connections $\nabla^\bullet,^\dagger\!\nabla^\bullet$ from the other.
\begin{lemma}\label{L:cond}
   Contravariant connections $\nabla^\bullet, ^\dagger\! \nabla^\bullet$ are associated if and only if in local coordinates their Christoffel symbols $\Gamma_k^{ij}$ and $^\dagger\Gamma_k^{ij}$, respectively, satisfy the following condition:
\begin{equation}\label{E:cond}
                    \pi^{ik} \Gamma_k^{jl} = \pi^{jk}\, ^\dagger\Gamma_k^{il}.
\end{equation}
\end{lemma}
\begin{proof}
It is enough to check (\ref{E:agree}) for $\alpha = dx^i$ and $\beta = dx^j$. Using (\ref{E:exact}) and the Jacobi identity for the Poisson tensor $\pi^{ij}$ we obtain that
\begin{eqnarray*}
 ^\dagger \nabla^{dx^i} \# dx^j - \nabla^{dx^j} \# dx^i - \#[dx^i,dx^j] = 
^\dagger\!\nabla^{dx^i} (\pi^{jk}\p_k) - \\
\nabla^{dx^j} (\pi^{ik}\p_k) - \frac{\p \pi^{ij}}{\p x^s}\pi^{sk} \p_k =
\pi^{is} \frac{\p \pi^{jk}}{\p x^s} \p_k + \pi^{jk}\, ^\dagger \Gamma_k^{il} \p_l -\\
\pi^{js} \frac{\p \pi^{ik}}{\p x^s} \p_k - \pi^{ik} \Gamma_k^{jl} \p_l - \frac{\p \pi^{ij}}{\p x^s} \pi^{sk} \p_k =  \left(\pi^{jk}\, ^\dagger \Gamma_k^{il} - \pi^{ik} \Gamma_k^{jl}\right)\p_l,
\end{eqnarray*}
whence the lemma follows.
\end{proof}
\begin{lemma}
A covariant connection $\nabla_\bullet$ respects the Poisson bivector field $\Pi$ if and only if the contravariant connection $^\dagger \nabla^\bullet$ induced by $^t \nabla_\bullet$ is torsion-free. 
\end{lemma}
\begin{proof}
 Denote by $T^{^\dagger\nabla}$ the torsion of the connection $^\dagger \nabla^\bullet$. Taking into account that $^\dagger \Gamma_k^{ij} = \pi^{is} \, ^t \Gamma^j_{sk} = \pi^{is} \Gamma^j_{ks}$ we get
\begin{align*}
   T^{^\dagger\nabla} (dx^i,dx^j) = \left(-^\dagger \Gamma_k^{ij}  + ^\dagger\! \Gamma_k^{ji}  - \frac{\p \pi^{ij}}{\p x^k}\right) dx^k = \\ -\left(\frac{\p \pi^{ij}}{\p x^k} + \Gamma^i_{ks} \pi^{sj} + \Gamma^j_{ks} \pi^{is}\right) dx^k = -\left(\nabla_{\p_k} \pi^{ij}\right)dx^k,
\end{align*}
whence the claim follows. 
\end{proof}
Thus a pair of associated contravariant connections $(\nabla^\bullet, ^\dagger\! \nabla^\bullet)$ such that $^\dagger\nabla^\bullet$ is torsion-free can be thought of as a natural counterpart of a covariant connection $\nabla_\bullet$ which respects the Poisson bivector field $\Pi$.
\begin{lemma}\label{L:resp}
 If contravariant connections $(\nabla^\bullet, ^\dagger\! \nabla^\bullet)$ are associated then for $\alpha,\beta,\gamma \in \Omega^1$
\begin{equation}\label{E:torfree}
\left(\nabla^\alpha \Pi\right)(\beta,\gamma) = - \Pi\left(\alpha, T^{^\dagger\nabla}(\beta,\gamma)\right).\end{equation}
In particular, if $^\dagger\nabla^\bullet$ is torsion-free then the contravariant connection $\nabla^\bullet$ respects the Poisson bivector field $\Pi$.
\end{lemma}
\begin{proof}
Using (\ref{E:cond}) we obtain (\ref{E:torfree}) by the following calculation:
\begin{align*}
    \nabla^{dx^m} \pi^{ij} = \pi^{mk} \frac{\p \pi^{ij}}{\p x^k} + \Gamma_k^{mi} \pi^{kj} + \Gamma_k^{mj} \pi^{ik} =\\
- \pi^{mk}\left(-^\dagger \Gamma_k^{ij}  + ^\dagger\! \Gamma_k^{ji}  - \frac{\p \pi^{ij}}{\p x^k}\right).
\end{align*}
\end{proof}

A contravariant connection $\nabla^\bullet$ which respects the Poisson bivector field $\Pi$ is called Poisson in \cite{RLF}. A connection $\nabla^\bullet$ is Poisson if and only if it satisfies the condition
\begin{equation}\label{E:basic}
                 \# \nabla^\alpha \beta = \nabla^\alpha \# \beta
\end{equation}
for any $\alpha, \beta \in \Omega^1(M)$. It was shown in \cite{RLF} that on any Poisson manifold there exist Poisson contravariant connections.

Assume once again that a pair $(\nabla^\bullet, ^\dagger\! \nabla^\bullet)$ of contravariant connections is induced by a pair $(\nabla_\bullet, ^t\!\nabla_\bullet)$ of covariant connections. If $\nabla_\bullet$ is torsion-free, then $\nabla^\bullet = ^\dagger\!\nabla^\bullet$. Thus a contravariant connection $\nabla^\bullet$ which satisfies the condition
\begin{equation}\label{E:self}
      \nabla^\alpha \# \beta - \nabla^\beta \# \alpha - \# [\alpha,\beta] =0        
\end{equation}
can be thought of as a natural counterpart of a torsion-free covariant connection.

If a contravariant connection $\nabla^\bullet$ satisfies condition (\ref{E:self}) and is torsion-free, it can be considered as an analogue of a torsion-free covariant connection which respects $\Pi$. 
\begin{lemma}\label{L:self}
  A torsion-free contravariant connection $\nabla^\bullet$ respects the Poisson bivector field $\Pi$ if and only if it satisfies (\ref{E:self}).
\end{lemma}
\begin{proof}
  Assume that $\nabla^\bullet$ respects $\Pi$ and therefore satisfies (\ref{E:basic}). Then (\ref{E:self}) immediately follows from (\ref{E:basic}) applied to the condition that $\nabla^\bullet$ is torsion-free. Conversely, if $\nabla^\bullet$ is torsion-free and satisfies (\ref{E:self}) it follows from Lemma \ref{L:resp} that $\nabla^\bullet$ respects $\Pi$. 
\end{proof}

It follows from Lemma \ref{L:self} that if $\nabla^\bullet$ is a torsion-free Poisson contravariant connection on $M$, then it is associated to itself and one can consider the pair $(\nabla^\bullet, \nabla^\bullet)$ as an example of a pair of associated connections $(\nabla^\bullet, ^\dagger\! \nabla^\bullet)$ such that $^\dagger \nabla^\bullet$ is torsion-free. Conversely, if contravariant connections $(\nabla^\bullet, ^\dagger\! \nabla^\bullet)$ on $M$ are associated and $^\dagger \nabla^\bullet$ is torsion-free then there exists a torsion-free Poisson connection on $M$.

\begin{proposition}
 If contravariant connections $(\nabla^\bullet, ^\dagger\! \nabla^\bullet)$ are associated and $^\dagger \nabla^\bullet$ is torsion-free then
\[
    \hat \nabla^\bullet = \frac{1}{3}\left(\nabla^\bullet + ^t\!\nabla^\bullet + ^\dagger\! \nabla^\bullet\right)
\]
is a torsion-free Poisson contravariant connection.
\end{proposition}
\begin{proof}
  It is clear that $\hat\nabla^\bullet$ is torsion-free. Taking into account Lemmas \ref{L:cond}, \ref{L:resp}, and \ref{L:self} we see that it remains to show that the Christoffel symbols of the connection $\hat \nabla^\bullet$,
\[
    \hat \Gamma_k^{ij} = \frac{1}{3}\left(\Gamma_k^{ij} + \Gamma_k^{ji} + \frac{\p \pi^{ji}}{\p x^k} + ^\dagger\!\Gamma_k^{ij}\right), 
\]
satisfy the condition
\[
    \pi^{ik} \hat\Gamma_k^{jl} = \pi^{jk} \hat \Gamma_k^{il},
\]
which immediately follows from (\ref{E:cond}), (\ref{E:torfree}), and the Jacobi identity for the Poisson tensor $\pi^{ij}$.
\end{proof}

On a symplectic manifold there exist symplectic (torsion-free) covariant connections. A symplectic covariant connection induces a torsion-free Poisson contravariant connection.

Another important example of torsion-free Poisson contravariant connections appears in the context of K\"ahler-Poisson manifolds. We call a complex manifold $M$ K\"ahler-Poisson if it is endowed with a Poisson tensor $g^{\bar lk}$ of type (1,1) with respect to the complex structure. Here $k$ and $\bar l$ are holomorphic and antiholomorphic indices, respectively, with respect to local holomorphic coordinates $\{z^k\}$ on $M$. The Jacobi identity for the Poisson tensor $g^{\bar lk}$ takes the form
\[
    g^{\bar ts}\frac{\p g^{\bar lk}}{\p z^s} = g^{\bar ls}\frac{\p g^{\bar tk}}{\p z^s} \mbox{ and } g^{\bar ts}\frac{\p g^{\bar lk}}{\p \bar z^t}=g^{\bar tk}\frac{\p g^{\bar ls}}{\p \bar z^t}.
\]
One can introduce a torsion-free Poisson contravariant connection $\nabla^\bullet$ on $M$ with the following Christoffel symbols:
\begin{equation}\label{E:christoff}
   \Gamma_m^{\bar lk} = - \frac{\p g^{\bar lk}}{\p z^m}, \ \Gamma_{\bar n}^{k\bar l} = \frac{\p g^{\bar lk}}{\p \bar z^n}
\end{equation}
(the symbols with the other types of indices are equal to zero). We call this contravariant connection K\"ahler-Poisson. If $g^{\bar lk}$ is nondegenerate, $M$ becomes a K\"ahler manifold and the contravariant connection $\nabla^\bullet$ is induced by the K\"ahler connection.

We are going to show that the Poisson manifold which is the dual of a semisimple Lie algebra does not admit torsion-free Poisson contravariant connections. First we need to prove a technical lemma.
\begin{lemma}\label{L:semisimple}
  On a semisimple Lie algebra $\mathfrak{g}$ any linear mapping $Q: \mathfrak{g} \to \mathfrak{g}$ such that 
\begin{equation}\label{E:symmetry}
   \left[ X, Q(Y)\right] =    \left[ Y, Q(X)\right]
\end{equation}
for all $X,Y \in \mathfrak{g}$ vanishes.
\end{lemma}
\begin{proof}
Let $Q: \mathfrak{g} \to \mathfrak{g}$ be a linear mapping satisfying (\ref{E:symmetry}). The Killing form $B(\cdot,\cdot)$ on $\mathfrak{g}$ is a nondegenerate invariant symmetric bilinear form. Thus $B([X,Y],Z) = B(X,[Y,Z])$, which implies that the function $\Theta(X,Y,Z) = B([X,Y], Q(Z)) = B(X,[Y,Q(Z)])$ is antisymmetric in $(X,Y)$ and symmetric in $(Y,Z)$. It can only happen if $\Theta = 0$. Since $B$ is nondegenerate and for $\mathfrak{g}$ semisimple $[\mathfrak{g},\mathfrak{g}] = \mathfrak{g}$, we obtain that $Q=0$.
\end{proof}

Let $\mathfrak{g}$ be a real finite-dimensional Lie algebra with a fixed basis $\{X^i\}$ and the structure constants $c_k^{ij}$ so that
\[
    [X^i,X^j] = c_k^{ij} X^k.
\]
Denote by $\{x^i\}$ the corresponding linear coordinates on the dual $\mathfrak{g}^*$ of the Lie algebra $\mathfrak{g}$.
The standard linear Poisson structure on $\mathfrak{g}^*$ is given by the formula
\begin{equation}\label{E:pic}
   \pi^{ij} = c_k^{ij} x^k.
\end{equation}

Assume that $\nabla^\bullet$ is a torsion-free Poisson contravariant connection on $\mathfrak{g}^*$ with the Christoffel symbols $\Gamma_k^{ij}$. Thus
\begin{equation}\label{E:gamma}
      - \Gamma_k^{ij} + \Gamma_k^{ji} - \frac{\p \pi^{ij}}{\p x^k} = 0 \mbox{ and } \pi^{ik} \Gamma_k^{jl} = \pi^{jk}\Gamma_k^{il}.
\end{equation}
Formulas (\ref{E:pic}) and (\ref{E:gamma}) imply that
\begin{equation}\label{E:zero}
      - \Gamma_k^{ij}(0) + \Gamma_k^{ji}(0) - c_k^{ij} = 0 \mbox{ and } c_s^{ik} \Gamma_k^{jl}(0) = c_s^{jk}\Gamma_k^{il}(0).
\end{equation}
Introduce a bilinear mapping $\Gamma: \mathfrak{g}\otimes\mathfrak{g} \to \mathfrak{g}$ such that
\[
 \Gamma(X^i,X^j) = - \Gamma_k^{ij}(0) X^k.
\]
It follows from (\ref{E:zero}) that the mapping $\Gamma(X,Y)$ satisfies the following conditions:
\begin{equation}\label{E:gammaxy}
   \Gamma(X,Y) - \Gamma(Y,X) -[X,Y]=0
\end{equation}
and
\begin{equation}\label{E:gammaxyz}
   [X, \Gamma(Y,Z)] = [Y, \Gamma(X,Z)].
\end{equation}

If $\mathfrak{g}$ is a semisimple Lie algebra it follows from Lemma \ref{L:semisimple} and formula (\ref{E:gammaxyz}) that for any element $Z\in \mathfrak{g}$ the mapping $\mathfrak{g} \ni X \to \Gamma(X,Z)$ vanishes, which implies that $\Gamma=0$. However, according to (\ref{E:gammaxy}) the skew-symmetric part of $\Gamma(X,Y)$ is $\frac{1}{2}[X,Y]$. This contradiction shows that on the dual of a semisimple Lie algebra endowed with the standard linear Poisson structure there are no torsion-free Poisson contravariant connections. It is interesting to compare this statement with the theorem by S. Gutt and J. Rawnsley that the dual of a semisimple Lie algebra does not admit a tangential deformation quantization (see \cite{GR} and also \cite{WTG}).

In the rest of the paper we will show that using Fedosov's approach one can give a geometric construction of a formal symplectic groupoid (see \cite{FSG}) on an arbitrary Poisson manifold endowed with a torsion-free Poisson connection. 

\medskip

\section{Nonlinear contravariant connections}

\medskip

The formal neighborhood $(X,Y)$ of a submanifold $Y$ of a manifold $X$ is the ringed space on $Y$ whose ring of global sections is $C^\infty(X)/\cap_{k=1}^\infty I_Y^k$, where $I_Y$ denotes the ideal of functions in $C^\infty(M)$ vanishing on $Y$ (see the Appendix to \cite{Deq}).
Let $p: E \to M$ be a vector bundle with a finite dimensional fibre. Denote by $\Vect (E,Z)$ the Lie algebra of formal vector fields on the formal neighborhood $(E,Z)$ of the zero section $Z$ of $E$. If $M$ is a Poisson manifold we call a mapping $D^\bullet: \Omega^1(M) \to \Vect(E,Z)$ a contravariant connection on the formal neighborhood $(E,Z)$ if it satisfies the axioms
\[
 D^{f\alpha} = p^*(f) D^\alpha, \quad [D^\alpha, p^* f] = p^*((\#\alpha)f),
\]
where $\alpha\in\Omega^1(M)$ and $f \in C^\infty(M)$. 

{\it Remark.} If $D^\bullet$ leaves invariant the space of fibrewise linear functions on $(E,Z)$, then it is induced by a linear connection on the dual bundle $E^* \to M$. In general, this is not the case and $D^\bullet$ will be referred to as a nonlinear contravariant connection on $(E,Z)$.

We will be interested in the nonlinear contravariant connections on the formal neighborhood $(\T,Z)$. 

If $\nabla^\bullet$ is a linear contravariant connection on $M$ with the Christoffel symbols $\Gamma_k^{ij}$, it induces a contravariant connection
$\bar\nabla^\bullet$ on $(\T,Z)$ expressed locally as
\[
    \bar\nabla^{dx^i} = \pi^{ij}\frac{\p}{\p x^j} + \Gamma_k^{ij} \xi_j \frac{\p}{\p \xi_k},
\]
where $\{x^i\}$ are local coordinates on $M$ lifted to $(T^*M,Z)$ and $\{\xi_i\}$ are the dual fibre coordinates.

For a global vector field $v\in\Vect(M)$ written locally as $v = v^i \p_i$ denote by $\hat v$ the multiplication operator by the global function on $C^\infty(T^*M,Z)$ whose local expression is $v^i(x)\xi_i$. For a global 1-form $\alpha \in \Omega^1(M)$ written locally as $\alpha = \alpha_i dx^i$ denote by $\hat \alpha$ the global operator on $C^\infty(T^*M,Z)$ written locally as $\alpha_i \frac{\p}{\p \xi_i}$.

Define the torsion $T^D$ of a nonlinear contravariant connection $D^\bullet$ on $(\T,Z)$ as the skew-symmetric bilinear mapping from $\Omega^1(M) \otimes_{C^\infty(M)} \Omega^1(M)$ to the space $\Vect_{\rm vert}(T^*M,Z)$ of vertical vector fields on the formal neighborhood $(T^*M,Z)$ given by the formula
\[
    T^D(\alpha,\beta) = \left[D^\alpha,\hat \beta\right] - \left[D^\beta,\hat \alpha\right] - \widehat{[\alpha,\beta]}.
\]
It is a morphism of $C^\infty(M)$-modules.
A nonlinear contravariant connection $D^\bullet$ on $(\T,Z)$ can be written in local coordinates as follows:
\begin{equation}\label{E:connloc}
   D^{dx^i} = \pi^{is}\frac{\p}{\p x^s} - A^i_s\frac{\p}{\p \xi_s},
\end{equation}
where $A^i_s = A^i_s(x,\xi)$ is a formal function (the fibre coordinates $\{\xi_i\}$ are treated as formal variables). Locally
\begin{equation}\label{E:torloc}
   T^D(dx^i,dx^j) = \frac{\p A^i_s}{\p \xi_j} \frac{\p}{\p \xi_s} - \frac{\p A^j_s}{\p \xi_i} \frac{\p}{\p \xi_s} -\frac{\p \pi^{ij}}{\p x^s} \frac{\p}{\p \xi_s}.
\end{equation}
One can check that the torsion $T^\nabla$ of a contravariant connection $\nabla^\bullet$ is related to the torsion $T^{\bar\nabla}$ of the induced contravariant connection $\bar\nabla^\bullet$ on $(\T,Z)$ as follows:
\[
     T^{\bar\nabla}(\alpha,\beta)  = \widehat{T^\nabla(\alpha,\beta)}
\]
for any $\alpha,\beta \in\Omega^1(M)$. 

A contravariant connection $\nabla^\bullet$ respects the Poisson bivector field $\Pi$ if and only if the induced mapping $\bar\nabla^\bullet$ satisfies the equation
\[
      \left[\left[\bar\nabla^\gamma, \widehat{\#\alpha}\right],\hat\beta\right] = \left[\left[\bar\nabla^\gamma, \hat\alpha\right],\widehat{\#\beta}\right]
\]
for any $\alpha,\beta,\gamma \in\Omega^1(M)$. We will call a nonlinear contravariant connection $D^\bullet$ Poisson if it satisfies the condition
\begin{equation}\label{E:dresppi} 
      \left[\left[D^\gamma, \widehat{\#\alpha}\right],\hat\beta\right] = \left[\left[D^\gamma, \hat\alpha\right],\widehat{\#\beta}\right]
\end{equation}
for any $\alpha,\beta,\gamma \in\Omega^1(M)$.

A simple calculation shows that $D^\bullet$ written in local coordinates as in (\ref{E:connloc}) satisfies (\ref{E:dresppi}) if and only if
\begin{equation}\label{E:dresppiloc}
    \pi^{ms} \frac{\p \pi^{ij}}{\p x^s} - \frac{\p A^m_s}{\p \xi_i} \pi^{sj} - \frac{\p A^m_s}{\p \xi_j} \pi^{is} =  0.
\end{equation}

Linear contravariant connections $(\nabla^\bullet,^\dagger\!\nabla^\bullet)$ are associated if and only if the mappings $\bar\nabla^\bullet, ^\dagger\!\bar\nabla^\bullet : \Omega^1(M) \to \Vect(E,Z)$ satisfy the condition
\begin{equation}\label{E:dagreeloc}
    \left[^\dagger\bar\nabla^\alpha, \widehat{\# \beta}\right] - \left[\bar\nabla^\beta, \widehat{\# \alpha}\right] - \widehat{\# [\alpha,\beta]} = 0.
\end{equation}
We will say that nonlinear contravariant connections $(D^\bullet,^\dagger\! D^\bullet)$ satisfying
\begin{equation}\label{E:dagree}
    \left[^\dagger D^\alpha, \widehat{\# \beta}\right] - \left[D^\beta, \widehat{\# \alpha}\right] - \widehat{\# [\alpha,\beta]} = 0
\end{equation}
for any $\alpha,\beta \in \Omega^1(M)$ are associated. One can check that if, locally, 
\begin{equation}\label{E:daggerloc}
D^{dx^i} = \pi^{is}\frac{\p}{\p x^s} - A^i_s\frac{\p}{\p \xi_s} \mbox{ and } ^\dagger D^{dx^i} = \pi^{is}\frac{\p}{\p x^s} - K^i_s\frac{\p}{\p \xi_s},
\end{equation}
then (\ref{E:dagree}) is equivalent to the condition
\begin{equation}\label{E:akagree}
     \pi^{is}A^j_s = \pi^{js} K^i_s.
\end{equation}
Using formulas (\ref{E:torloc}), (\ref{E:dresppiloc}), and (\ref{E:akagree}) one can prove the following lemma. 
\begin{lemma}\label{L:nonlinresp}
 If nonlinear contravariant connections $(D^\bullet,^\dagger\! D^\bullet)$ are associated and $^\dagger D^\bullet$ is torsion-free, then the nonlinear contravariant connection $D^\bullet$ is Poisson.
\end{lemma}

The Poisson structure $\Pi$ on $M$ induces a fibrewise presymplectic form $\Omega$ on $T^*M$ given locally by the formula
\[
           \Omega = \frac{1}{2} \pi^{ij}d\xi_i \wedge d\xi_j.
\]
Denote by $\H_\Omega$ the subspace of $C^\infty(T^*M,Z)$ of Hamiltonian functions of the fibrewise presymplectic form $\Omega$. Namely, $F \in \H_\Omega$ if there is a vertical formal vector field $H_F\in\Vect_{\rm vert}(T^*M,Z)$ such that
\[
     \iota(H_F)\Omega = - d_{\xi} F,
\]
where $d_\xi$ is the fibrewise differential. The choice of $H_F$ is, in general, not unique.
If locally $H_F = a_i \frac{\p}{\p \xi_i}$, then
\begin{equation}\label{E:ham}
         \frac{\p F}{\p \xi_i} = - \pi^{ji} a_j = \pi^{ij} a_j. 
\end{equation}
It follows from (\ref{E:ham}) that
\begin{equation}\label{E:hameqn}
               \pi^{ij}\frac{\p a_j}{\p \xi_k} = \pi^{kj}\frac{\p a_j}{\p \xi_i}.
\end{equation}
The space $\H_\Omega$ is a Poisson algebra with respect to the pointwise product and the bracket $\{\cdot,\cdot\}_\Omega$ defined as follows:
\begin{equation}\label{E:bracketomega}
   \{F,G\}_\Omega = H_F G = - H_G F =  a_i \frac{\p G}{\p \xi_i} = \pi^{ij}a_i b_j,            
\end{equation}
where $F,G \in \H_\Omega$ and, locally, $H_F = a_i \frac{\p}{\p \xi_i}$ and $H_G = b_i \frac{\p}{\p \xi_i}$.

{\it Remark.} For any formal function $F\in C^\infty(TM,Z)$ its pullback $F^\#$ with respect to the mapping $\#$ belongs to $\H_\Omega$. For $F=F(x,y)$ we have $F^\#(x,\xi) = F(x,\xi_i\pi^{i\cdot})$. Therefore
\[
    \frac{\p F^\#}{\p \xi_i} = \pi^{ij}\left(\frac{\p F}{\p y^j}\right)^\#,
\]
where $\{y^i\}$ are the local fibre coordinates on $TM$ dual to $\{\xi_i\}$.

\begin{proposition}\label{P:dderiv}
  Given a Poisson nonlinear contravariant connection $D^\bullet$, the formal vector fields  $D^\alpha, \alpha \in \Omega^1(M),$ are derivations of the Poisson algebra $(\H_\Omega, \{\cdot,\cdot\}_\Omega)$.
\end{proposition}

\begin{proof}
 For $F \in \H_\Omega$ assume that locally $H_F = a_i \frac{\p}{\p \xi_i}$ and formula (\ref{E:ham}) holds. In order to prove the Proposition we have to show that
\[
    b_t \frac{\p}{\p \xi_t} := [D^{dx^i}, H_F],
\]
is a Hamiltonian vector field of the local Hamiltonian function $D^{dx^i}F$, i.e., that
\begin{equation}\label{E:locpotentbt}
    \frac{\p}{\p \xi_j} \left(D^{dx^i}F\right) = \pi^{jt} b_t.
\end{equation}
We have, using (\ref{E:ham}),
\begin{align}\label{E:one}
   \frac{\p}{\p \xi_j}\left(D^{dx^i}F\right) = \frac{\p}{\p \xi_j} \left(\pi^{is}\frac{\p F}{\p x^s} - A^i_s \frac{\p F}{\p \xi_s}\right) = \pi^{is} \frac{\p^2 F}{\p x^s \p \xi_j} -\frac{\p A^i_s}{\p \xi_j} \frac{\p F}{\p \xi_s} - \nonumber\\
 A^i_s \frac{\p^2 F}{\p \xi_s \p \xi_j} = \pi^{is} \frac{\p}{\p x^s}\left(\pi^{jt} a_t\right) -\frac{\p A^i_s}{\p \xi_j} \pi^{st} a_t - A^i_s \pi^{jt} \frac{\p a_t}{\p \xi_s} = \\
\pi^{is} \pi^{jt} \frac{\p a_t}{\p x^s} + \pi^{is} \frac{\p \pi^{jt}}{\p x^s} a_t -\frac{\p A^i_s}{\p \xi_j} \pi^{st} a_t  - A^i_s \pi^{jt} \frac{\p a_t}{\p \xi_s}.\nonumber
\end{align}
On the other hand we obtain that
\begin{align}\label{E:two}
  b_t \frac{\p}{\p \xi_t} = \left[\pi^{is}\frac{\p }{\p x^s} - A^i_s\frac{\p }{\p \xi_s}, a_t \frac{\p}{\p \xi_t}\right] = \\
\pi^{is} \frac{\p a_t}{\p x^s}\frac{\p}{\p \xi_t}  - A^i_s \frac{\p a_t}{\p \xi_s}\frac{\p}{\p \xi_t} + \frac{\p A^i_t}{\p \xi_k} a_k \frac{\p}{\p \xi_t}.\nonumber
\end{align}
Taking into account that $D^\bullet$ is a Poisson connection we obtain formula (\ref{E:locpotentbt}) from (\ref{E:dresppiloc}), (\ref{E:one}), and (\ref{E:two}), whence the Proposition follows.
\end{proof}

\section{Lifting functions via formal nonlinear connections}

There is a natural grading $\deg$ on the space $C^\infty(T^*M,Z)$ by the powers of the formal fibre variables $\xi_i$. It induces a grading on the endomorphisms of $C^\infty(T^*M,Z)$. Thus $\deg(\xi_i) = 1$ and $\deg\left(\frac{\p}{\p \xi_i}\right) = -1$.
For a linear contravariant connection $\nabla^\bullet$ on $M$ and a 1-form $\alpha \in \Omega^1(M)$ we have $\deg\left(\bar\nabla^\alpha\right) =0$. A nonlinear contravariant connection $D^\bullet$ can be expanded into a $\xi$-adically convergent series 
\[
   D^\bullet = D^\bullet_{-1} + D^\bullet_0 + \ldots,
\]
where $D^\bullet_d$ is the homogeneous component of $D^\bullet$ of degree $d$. There is a unique contravariant connection $\nabla^\bullet$ such that $D^\bullet_0 = \bar\nabla^\bullet$. For any $d \neq 0$ the component $D^\bullet_d$ is a mapping from $\Omega^1(M)$ to $\Vect_{\rm vert}(T^*M,Z)$. 
The mapping $\alpha \mapsto D^\alpha_{-1}$ determines uniquely a linear endomorphism $\psi_D$ of the cotangent bundle $T^*M$ such that 
\[
   D_{-1}^\alpha = \widehat{\psi_D(\alpha)}           
\]
for all $\alpha \in \Omega^1(M)$. We will call a nonlinear contravariant connection $D^\bullet$ {\it invertible} if $\psi_D$ is invertible. We will give necessary and sufficient conditions on an invertible nonlinear contravariant connection $D^\bullet$ under which the system 
\begin{equation}\label{E:lift}
             D^\bullet F = 0
\end{equation}
with the initial condition $F|_Z = p^*(f)|_Z$ has a unique solution $F \in C^\infty(T^*M,Z)$
for an arbitrary $f \in C^\infty(M)$.

The curvature of a nonlinear contravariant connection $D^\bullet$ is a skew-symmetric morphism of $C^\infty(M)$-modules
\[
   R^D: \Omega^1(M) \otimes_{C^\infty(M)} \Omega^1(M) \to \Vect_{\rm vert}(T^*M,Z)
\] 
such that for $\alpha,\beta \in \Omega^1(M)$
\[
    R^D(\alpha,\beta) = \left[D^\alpha,D^\beta\right] - D^{[\alpha,\beta]}.
\]
If a nonlinear contravariant connection $D^\bullet$ is written locally as in (\ref{E:connloc}), then 
\begin{equation}\label{E:curvloc1}
   R^D(dx^i,dx^j) =\left[D^{dx^i},D^{dx^j}\right] - \frac{\p \pi^{ij}}{\p x^k} D^{dx^k}= R^{ij}_k \frac{\p}{\p \xi_k},
\end{equation}
where
\begin{equation}\label{E:curvloc2}
   R^{ij}_k = -D^{dx^i}\left(A^j_k\right) + D^{dx^j}\left(A^i_k\right)
+ \frac{\p \pi^{ij}}{\p x^s} A^s_k.
\end{equation}
System (\ref{E:lift}) written locally takes the form
\begin{equation}\label{E:liftloc}
              \pi^{is}\frac{\p F}{\p x^s} = A^i_s \frac{\p F}{\p \xi_s}.
\end{equation}
Since $D^\bullet$ is invertible, the matrix $\left(A^i_j\right)$ is formally invertible.
Denote by $\left(B^i_j\right)$ its inverse. System (\ref{E:liftloc}) is equivalent to the following one:
\begin{equation}\label{E:bpiloc}
           \frac{\p F}{\p \xi_i} = B^i_j \pi^{js} \frac{\p F}{\p x^s}.
\end{equation}
It is well known that system (\ref{E:bpiloc}) has a unique formal solution $F(x,\xi)$ with the initial condition $F(x,0) = f(x)$ for an arbitrary smooth function $f(x)$ if and only if the operators
\[
    B^i_j D^{dx^j} =  B^i_j \pi^{js} \frac{\p }{\p x^s} - \frac{\p }{\p \xi_i}
\]
pairwise commute. We have, taking into account the identity 
\[
D^{dx^k}B^j_l = - B^j_p D^{dx^k}\left(A^p_q\right) B^q_l
\]
and formulas (\ref{E:curvloc1}) and (\ref{E:curvloc2}), that
\begin{align*}
   \left[B^i_k D^{dx^k},B^j_l D^{dx^l}\right] = B^i_k D^{dx^k}\left(B^j_l\right)D^{dx^l} - 
B^j_k D^{dx^k}\left(B^i_l\right)D^{dx^l}+\\ 
B^i_k B^j_l\left[D^{dx^k},D^{dx^l}\right] = 
B^i_k B^j_l\Big(- D^{dx^k}\left(A^l_q\right)B^q_pD^{dx^p} + D^{dx^l}\left(A^k_q\right)B^q_pD^{dx^p} + \\
\left[D^{dx^k},D^{dx^l}\right]\Big) = B^i_k B^j_l\left(\left(R^{kl}_q - \frac{\p \pi^{kl}}{\p x^s}A^s_q\right)B^q_p D^{dx^p} + \left[D^{dx^k},D^{dx^l}\right]\right) = \\
B^i_k B^j_l\left(R^{kl}_q B^q_p D^{dx^p} + R^{kl}_q \frac{\p}{\p \xi_q}\right) = B^i_k B^j_l R^{kl}_q B^q_p \pi^{ps} \frac{\p}{\p x^s}.
\end{align*}
Since the matrix $\left(B^i_j\right)$ is invertible, we obtain the following proposition.
\begin{proposition}\label{P:system}
  For an invertible nonlinear contravariant connection $D^\bullet$ written locally as (\ref{E:connloc}) system (\ref{E:lift}) with the initial condition $F(x,0) = f(x)$ has a unique solution $F = F(x,\xi)$ for an arbitrary smooth function $f(x)$ if and only if \begin{equation}\label{E:loccrit}
     R^{kl}_q B^q_p \pi^{ps} =0,
\end{equation}
where $R^{kl}_q$ are the coefficients of the curvature of the connection $D^\bullet$ given by formula (\ref{E:curvloc2}) and $\left(B^i_j\right)$ is the inverse of $\left(A^i_j\right)$.
\end{proposition}  
Proposition \ref{P:system} gives only a local necessary and sufficient condition of 
the existence and uniqueness of solutions of system (\ref{E:lift}). To give a global criterion, we assume additionally that the invertible connection $D^\bullet$ has an associated invertible nonlinear contravariant connection $^\dagger D^\bullet$.
Suppose that $D^\bullet$ and $^\dagger D^\bullet$ are written locally as in (\ref{E:daggerloc}). Since $D^\bullet$ and $^\dagger D^\bullet$ are invertible, there exist the inverse matrices $(B^i_j)$ and $(L^i_j)$ of $(A^i_j)$ and $(K^i_j)$, respectively. It follows from (\ref{E:akagree}) that
\begin{equation}\label{E:blagree}
     \pi^{is}B^j_s = \pi^{js} L^i_s.
\end{equation}
Then, if $R^{kl}_q$ are the coefficients of the curvature $R^D$, condition (\ref{E:loccrit})  is equivalent to the following one:
\begin{equation}\label{E:intercrit}
     R^{kl}_q \pi^{qs} L^p_s =0.
\end{equation}
Since $(L^i_j)$ is invertible, we obtain that (\ref{E:intercrit}) is equivalent to the condition
\[
   R^{kl}_q \pi^{qs} =0
\]
which means that the vertical formal vector field $R^D(dx^k,dx^l)$ annihilates the fibrewise presymplectic form $\Omega$. Thus we arrive at the following theorem.
\begin{theorem}\label{T:system}
  Let $(D^\bullet,^\dagger D^\bullet)$  be a pair of associated invertible nonlinear formal contravariant connections. Then system (\ref{E:lift}) with the initial condition $F|_Z = p^*(f)|_Z$ has a unique solution $F \in C^\infty(T^*M,Z)$ for an arbitrary $f \in C^\infty(M)$ if and only if for any $\alpha,\beta \in \Omega^1(M)$ the vertical vector field $R^D(\alpha,\beta)$ is in the kernel of the fibrewise presymplectic form $\Omega$,
\[
    \iota\left(R^D(\alpha,\beta)\right) \Omega = 0.
\]
\end{theorem}

Another important property of such a pair of contravariant connections is given by the following lemma.

\begin{lemma}\label{L:kernel}
   Let $(D^\bullet,^\dagger D^\bullet)$  be a pair of associated invertible nonlinear formal contravariant connections. Then the kernel of the connection $D^\bullet$ is a subspace of $\H_\Omega$.
\end{lemma}
\begin{proof}
 
Write the connections $D^\bullet$ and $^\dagger D^\bullet$ locally as in (\ref{E:daggerloc}) and let $(B^i_j)$ and $(L^i_j)$ be the inverse matrices of $(A^i_j)$ and $(K^i_j)$, respectively. If $F\in C^\infty(T^*M,Z)$ is a solution of system (\ref{E:lift}), then $F$ is also a solution of system  (\ref{E:bpiloc}). It follows form formula (\ref{E:blagree}) that 
\[
           \frac{\p F}{\p \xi_i} = - \pi^{ij} L^s_j\frac{\p F}{\p x^s},        
\]
which implies that $F \in \H_\Omega$.
\end{proof}

In general, system (\ref{E:lift}) can have the same space of solutions for different nonlinear contravariant connections $D^\bullet$. Given two pairs of associated invertible nonlinear contravariant connections, $(D^\bullet_1, {^\dagger\!} D^\bullet_1)$ and  $(D^\bullet_2, {^\dagger\!} D^\bullet_2)$, one can show that if for any $\alpha,\beta \in \Omega^1(M)$
\[
     \left[D_1^\alpha, \widehat{\# \beta}\right] = \left[D_2^\alpha, \widehat{\# \beta}\right],
\]
then the kernels of the connections $D^\bullet_1, D^\bullet_2$ coincide.

\section{The Fundamental Equation}

In this section we will assume that $(\nabla^\bullet, ^\dagger\! \nabla^\bullet)$ 
are associated global contravariant connections with $^\dagger \nabla^\bullet$ torsion-free on a Poisson manifold $(M,\Pi)$. Starting with these connections, we are going to construct a global Poisson morphism $\theta: (C^\infty(M), \{\cdot,\cdot\}) \to (\H_\Omega, \{\cdot,\cdot\}_\Omega)$.

Consider local nonlinear contravariant connections
\begin{align}\label{E:ncd}
   D^{dx^i} = \pi^{is}\frac{\p}{\p x^s} + \Gamma^{ij}_s \xi_j \frac{\p}{\p \xi_s} - v^i_s \frac{\p}{\p \xi_s}, \\
^\dagger\! D^{dx^i} = \pi^{is}\frac{\p}{\p x^s} + ^\dagger\! \Gamma^{ij}_s \xi_j \frac{\p}{\p \xi_s} - w^i_s \frac{\p}{\p \xi_s},
\end{align}
where $\Gamma^{ij}_k,^\dagger\! \Gamma^{ij}_k$ are the Christoffel symbols of the connections $\nabla^\bullet,^\dagger\!\nabla^\bullet$, respectively, and 
$v_k^i = v_k^i(x,\xi), w_k^i = w_k^i(x,\xi)$ are formal functions.

It follows from formulas (\ref{E:dagreeloc}) and (\ref{E:dagree}) that connections $D^\bullet$ and $^\dagger D^\bullet$ are associated if and only if
\[
     \pi^{is}v_s^j = \pi^{js}w_s^i.
\]
Connection $^\dagger D^\bullet$ is torsion-free if and only if
\[
     \frac{\p w_k^i}{\p \xi_j} - \frac{\p w_k^j}{\p \xi_i} = 0,
\]
or, equivalently, if there exists a formal function $u_k(x,\xi)$ such that
\[
    w_k^i = \frac{\p u_k}{\p \xi_i}.
\]
Now assume that $D^\bullet$ and $^\dagger D^\bullet$ are associated and $^\dagger D^\bullet$ is torsion-free. We will write
\begin{equation}\label{E:ncdb}
  ^\dagger\! D^{dx^i} = \pi^{is}\frac{\p}{\p x^s} + ^\dagger\! \Gamma^{ij}_s \xi_j \frac{\p}{\p \xi_s} - \frac{\p u_s}{\p \xi_i} \frac{\p}{\p \xi_s}.
\end{equation}
The condition that $D^\bullet$ and $^\dagger D^\bullet$ are associated takes the form
\begin{equation}\label{E:vuagree}
     \pi^{is}v_s^j = \pi^{js}\frac{\p u_s}{\p \xi_i} = \frac{\p \left(\pi^{js} u_s\right)}{\p \xi_i}.     
\end{equation}
It follows from (\ref{E:vuagree}) that $\pi^{js} u_s$ is a local Hamiltonian function for the fibrewise presymplectic structure defined by $\Omega$.

The coefficients $\bar R^{ij}_k$ of the curvature $R^{\bar \nabla}$ of the mapping  $\bar \nabla^\bullet$ are given by the formula
\begin{equation}\label{E:barr}
     \bar R^{ij}_k = \pi^{is}\frac{\p \Gamma^{jp}_k}{\p x^s}\xi_p - \pi^{js}\frac{\p \Gamma^{ip}_k}{\p x^s}\xi_p - \frac{\p \pi^{ij}}{\p x^s} \Gamma^{sp}_k \xi_p - \Gamma^{iq}_k \Gamma^{jp}_q \xi_p +  \Gamma^{jq}_k \Gamma^{ip}_q \xi_p.
\end{equation}
Denote by $\E$ the fibrewise Euler vector field on $T^*M$ expressed locally as
\[
     \E = \xi_k\frac{\p}{\p \xi_k}
\]
and introduce a skew-symmetric bilinear morphism of $C^\infty(M)$-modules 
\[
   \bar Q : \Omega^1(M) \otimes_{C^\infty(M)} \Omega^1(M) \to C^\infty(\T)
\]
as follows:
\[
   \bar Q (\alpha,\beta) = \frac{1}{2} 
\Omega\left(\E, R^{\bar\nabla}(\alpha,\beta)\right).
\]
Introduce a function $\bar Q^{ij} = \bar Q(dx^i,dx^j)$. Then
\begin{equation}\label{E:barqijdef}
     \bar Q^{ij} = \frac{1}{2} \xi_t \pi^{tk} \bar R^{ij}_k.
\end{equation}
Using formulas (\ref{E:cond}), (\ref{E:torfree}), (\ref{E:barqijdef}), and the Jacobi identity for the Poisson tensor $\pi^{ij}$ one can show that
\begin{align}\label{E:gggg}
    \frac{\p}{\p \xi_p}\left(\pi^{tk}\bar R^{ij}_k\right) = \frac{1}{2}\pi^{is}\pi^{jk}\left( \frac{\p}{\p x^s} \left(^\dagger \Gamma_k^{tp} + ^\dagger\! \Gamma_k^{pt}\right) - \frac{\p}{\p x^k} \left(^\dagger \Gamma_s^{tp} + ^\dagger\! \Gamma_s^{pt}\right)\right) + \nonumber\\
\Gamma^{it}_k \pi^{kl} \Gamma^{jp}_l + \Gamma^{ip}_k \pi^{kl} \Gamma^{jt}_l.
\end{align}
The right-hand side of (\ref{E:gggg}) is symmetric with respect to the permutation of the indices $p$ and $t$. It follows from (\ref{E:barqijdef}) and (\ref{E:gggg}) that
\begin{equation}\label{E:barq}
    \pi^{tk}\bar R^{ij}_k = \frac{\p}{\p \xi_t} \bar Q^{ij}.
\end{equation}
Formula (\ref{E:barq}) means that for any indices $i,j$ the function $\bar Q^{ij}$ is a local Hamiltonian function for the presymplectic structure defined by $\Omega$. This can be written globally as follows:
\[
     \iota\left(R^{\bar\nabla}(\alpha,\beta)\right)\Omega = - d_{\xi} \bar Q(\alpha,\beta).
\]
The coefficients $R^{ij}_k$ of the curvature $R^D$ of the connection $D^\bullet$ are expressed as
\begin{align}\label{E:rijk}
   R^{ij}_k =  \bar R^{ij}_k - \pi^{is}\frac{\p v^j_k}{\p x^s} + \pi^{js}\frac{\p v^i_k}{\p x^s} + \frac{\p \pi^{ij}}{\p x^s} v^s_k - \Gamma^{ip}_s \xi_p \frac{\p v^j_k}{\p \xi_s} + \nonumber\\
\Gamma^{jp}_s \xi_p \frac{\p v^i_k}{\p \xi_s} + \Gamma^{ip}_k v^j_p - \Gamma^{jp}_k v^i_p + v^i_s \frac{\p v^j_k}{\p \xi_s} - v^j_s \frac{\p v^i_k}{\p \xi_s}.
\end{align}
We want to show that there is a  function $Q^{ij}$ such that
\begin{equation}\label{E:q}
    \pi^{tk}R^{ij}_k = \frac{\p}{\p \xi_t} Q^{ij}.
\end{equation}
Taking into account Lemma \ref{L:resp} we get from (\ref{E:rijk}) that
\begin{align}\label{E:pirijk}
    \pi^{tk}R^{ij}_k = \pi^{tk}\bar R^{ij}_k - \pi^{is}\frac{\p \left(\pi^{tk}v^j_k\right)}{\p x^s} + \pi^{js}\frac{\p \left(\pi^{tk}v^i_k\right)}{\p x^s} + \frac{\p \pi^{ij}}{\p x^s} \pi^{tk}v^s_k - \nonumber\\
\Gamma^{ip}_s \xi_p \frac{\p \left(\pi^{tk}v^j_k\right)}{\p \xi_s} + 
\Gamma^{jp}_s \xi_p \frac{\p \left(\pi^{tk}v^i_k\right)}{\p \xi_s} -
\Gamma^{it}_s \pi^{sk}v^j_k + \Gamma^{jt}_s \pi^{sk}v^i_k + \\
v^i_s \pi^{tk}\frac{\p v^j_k}{\p \xi_s} - v^j_s \pi^{tk}\frac{\p v^i_k}{\p \xi_s}. \nonumber
\end{align}
Using formulas (\ref{E:hameqn}), (\ref{E:vuagree}), and (\ref{E:barq}) we obtain from (\ref{E:pirijk}) that
\begin{align}\label{E:vtou}
    \pi^{tk}R^{ij}_k = \frac{\p}{\p \xi_t}\Big(\bar Q^{ij}- \pi^{is}\frac{\p \left(\pi^{jk}u_k\right)}{\p x^s} + \pi^{js}\frac{\p \left(\pi^{ik}u_k\right)}{\p x^s} + \frac{\p \pi^{ij}}{\p x^s} \pi^{sk}u_k - \nonumber\\
\Gamma^{ip}_s \xi_p \frac{\p \left(\pi^{jk}u_k\right)}{\p \xi_s} + \Gamma^{jp}_s \xi_p \frac{\p \left(\pi^{ik}u_k\right)}{\p \xi_s} + \pi^{st} v^i_s v^j_t  \Big).
\end{align}
Formulas (\ref{E:bracketomega}) and  (\ref{E:vtou}) imply that (\ref{E:q}) holds for
\begin{align}\label{E:qij}
    Q^{ij} = \bar Q^{ij}- 
\bar\nabla^i \left(\pi^{jk}u_k\right) + \bar\nabla^j \left(\pi^{ik}u_k\right) + \nonumber\\
\frac{\p \pi^{ij}}{\p x^s} \pi^{sk}u_k +\{\pi^{ik} u_k,\pi^{jl} u_l\}_\Omega.
\end{align}

Let $U$ be a global vertical 1-form on $(\T,Z)$. For a 1-form $\alpha \in \Omega^1(M)$ define an element $U(\alpha) \in C^\infty(\T,Z)$ as follows:
\[
            U(\alpha) = \langle U, \hat \alpha \rangle.   
\]
The mapping $\Omega^1(M) \ni \alpha \mapsto U(\alpha)$ from $\Omega^1(M)$ to $C^\infty(\T,Z)$ is a global morphism of $C^\infty(M)$-modules.
Now assume that $U$ is a global vertical 1-form on $(\T,Z)$ for which there exist local potentials $u_k(x,\xi)$ and $v^i_k(x,\xi)$ (on a covering by coordinate charts) satisfying (\ref{E:vuagree}) and such that 
\begin{equation}\label{E:udxi}
     U(dx^i) = \pi^{ik}u_k, 
\end{equation}
\begin{equation}\label{E:invert}
\ u_k = -\xi_k \pmod{\xi^2}, \mbox{ and } v^i_k = \delta^i_k \pmod{\xi}.
\end{equation}
We see from (\ref{E:vuagree}) and (\ref{E:udxi}) that $U(\alpha) \in \H_\Omega$ for any $\alpha \in \Omega^1(M)$.

Introduce {\it local} contravariant connections $D^\bullet$ and $^\dagger D^\bullet$ by formulas (\ref{E:ncd}) and (\ref{E:ncdb}), respectively. Thus $^\dagger D^\bullet$ is torsion-free. It follows from (\ref{E:vuagree}) that $D^\bullet$ and $^\dagger D^\bullet$ are associated. Conditions (\ref{E:invert}) imply that $D^\bullet$ and $^\dagger D^\bullet$ are both invertible. Lemma \ref{L:kernel} means that system (\ref{E:lift}) with $F\in C^\infty(\T,Z)$ has the same solutions as system (\ref{E:lift}) with $F \in \H_\Omega$. It follows from (\ref{E:bracketomega}) and (\ref{E:vuagree}) that if $F \in \H_\Omega$, system (\ref{E:lift}) can be rewritten as
\begin{equation}\label{E:hamilt}
    \bar\nabla^\alpha F = \{U(\alpha),F\}_\Omega.
\end{equation}
Introduce a global skew-symmetric mapping
\[
   Q : \Omega^1(M) \otimes_{C^\infty(M)} \Omega^1(M) \to C^\infty(\T,Z)
\]
by the formula
\begin{align}\label{E:qijglob}
    Q(\alpha,\beta) = \bar  Q(\alpha,\beta) - \nabla^\alpha U(\beta) + \nabla^\beta U(\alpha) + \nonumber\\
  U([\alpha,\beta]) + \{U(\alpha), U(\beta)\}_\Omega.
\end{align}
It is easy to check that $Q$ is a morphism of $C^\infty(M)$-modules. Notice that
\[
      Q(dx^i,dx^j) = Q^{ij},
\] 
where $Q^{ij}$ is given by formula (\ref{E:qij}). It follows from formula (\ref{E:q}) that
\begin{equation}\label{E:rdq}
     \iota\left(R^D(\alpha,\beta)\right)\Omega = - d_{\xi} Q(\alpha,\beta).
\end{equation}
We obtain from Theorem \ref{T:system} and (\ref{E:rdq}) that for any $f \in C^\infty(M)$
system (\ref{E:hamilt}) has a unique solution $F \in \H_\Omega$ such that $F|_Z = p^*(f)$ if and only if
\begin{equation}\label{E:qzero}
   d_{\xi} Q(\alpha,\beta) = 0
\end{equation}
Taking into account conditions (\ref{E:invert}) we get from (\ref{E:qzero}) that
\[
    Q(\alpha,\beta) = \Pi(\alpha,\beta).
\]
Summing up the results obtained in this section we can state the following theorem.
\begin{theorem}\label{T:lift}
If $U$ is a global vertical 1-form on $(\T,Z)$ which has local potentials $u_k(x,\xi)$ and $v^i_k(x,\xi)$ such that conditions (\ref{E:vuagree}), (\ref{E:udxi}), and (\ref{E:invert}) hold and which satisfies the equation
\begin{align}\label{E:Fundamental}
    \Pi(\alpha,\beta) = \bar  Q(\alpha,\beta) - \bar\nabla^\alpha U(\beta) + \bar\nabla^\beta U(\alpha) + \nonumber\\
  U([\alpha,\beta]) + \{U(\alpha), U(\beta)\}_\Omega,
\end{align}
then for any $f \in C^\infty(M)$ system (\ref{E:hamilt}) has a unique solution $F \in \H_\Omega$ such that $F|_Z = p^*(f)$.
\end{theorem}
We call equation (\ref{E:Fundamental}) the Fundamental Equation.
If $U$ is as in Theorem \ref{T:lift}, introduce the mapping 
\[
\theta: C^\infty(M) \to \H_\Omega
\]
which maps $f \in C^\infty(M)$ to the solution $F\in \H_\Omega$ of system (\ref{E:hamilt}) with the initial condition $F|_Z = p^*(f)$. We get from Lemma \ref{L:nonlinresp} and Proposition \ref{P:dderiv} that the mapping
\[
   \H_\Omega \ni F \mapsto \bar\nabla^\alpha F - \{U(\alpha),F\}_\Omega
\]
is a derivation of the Poisson algebra $\left(\H_\Omega,\{\cdot,\cdot\}_\Omega\right)$, whence we see that the solutions of system (\ref{E:hamilt}) are closed with respect to the Poisson bracket $\{\cdot,\cdot\}_\Omega$. Using this fact one can readily show that the mapping $\theta$ is a Poisson morphism.
\begin{proposition}\label{P:poisson}
 If $U$ is a global vertical 1-form on $(\T,Z)$ which has local potentials $u_k(x,\xi)$ and $v^i_k(x,\xi)$ such that conditions (\ref{E:vuagree}), (\ref{E:udxi}), and (\ref{E:invert}) hold and which satisfies the Fundamental Equation, the corresponding mapping $\theta$ is a Poisson morphism from $(C^\infty(M),\{\cdot,\cdot\})$ to $\left(\H_\Omega,\{\cdot,\cdot\}_\Omega\right)$.
\end{proposition}

The following theorem which states the existence and uniqueness of a normalized solution of the Fundamental Equation can be proved by induction. 
\begin{theorem}\label{T:Fundamentalsolution}
  Given a pair $(\nabla^\bullet, ^\dagger\!\nabla^\bullet)$ of associated contravariant connections  on $M$ such that $^\dagger\!\nabla^\bullet$ is torsion-free,
there exists a unique global solution $U$ of the Fundamental Equation (\ref{E:Fundamental}) which satisfies the normalization condition
\begin{equation}\label{E:normcond}
              \langle U, \E \rangle = 0
\end{equation}
and has local potentials $u_k(x,\xi)$ and $v^i_k(x,\xi)$ such that conditions (\ref{E:vuagree}), (\ref{E:udxi}), and (\ref{E:invert}) hold.
\end{theorem}

We will give a recurrence formula for 
\begin{equation}\label{E:solutionu}
           u^i := U(dx^i)
\end{equation}
which proves the uniqueness part of the theorem.
For a function $\alpha(x,\xi) \in C^\infty(\T,Z)$ denote by $\alpha^{(s)}$ its homogeneous component of degree $s$ with respect to the fibre variables $\xi$. Assume that $U$ is a normalized solution of the Fundamental Equation. Locally normalization condition (\ref{E:normcond}) reads
\begin{equation}\label{E:normcondloc}
   u^i\xi_i = 0,
\end{equation}
where $u^i$ is given by (\ref{E:solutionu}).
The Fundamental Equation can be written locally as follows:
\begin{equation}\label{E:Fundamentalloc}
   \pi^{ij} = \bar  Q^{ij} - \bar\nabla^{dx^i} u^j + \bar\nabla^{dx^j} u^i + \frac{\p \pi^{ij}}{\p x^k} u^k + \{u^i, u^j\}_\Omega.
\end{equation}
Taking into account formulas (\ref{E:bracketomega}), (\ref{E:barqijdef}), and (\ref{E:invert}), observe that $\bar  Q^{ij}$ is quadratic in $\xi$, the Poisson bracket $\{\cdot,\cdot\}_\Omega$ reduces the degree of homogeneity by two, $(u^i)^{(1)} = -\pi^{ij}\xi_j$, and for $F \in \H_\Omega$
\begin{equation}\label{E:firstterm}
    \{(u^i)^{(1)}, F\}_\Omega = \{-\pi^{ij}\xi_j, F\}_\Omega = \frac{\p F}{\p \xi_i}.
\end{equation}
Extracting from (\ref{E:Fundamentalloc}) the homogeneous component of degree $s-1, s\geq 2,$ we obtain, using (\ref{E:firstterm}), the following formula:
\begin{equation}\label{E:uniqueu}
   \frac{\p (u^i)^{(s)}}{\p \xi_j} - \frac{\p (u^j)^{(s)}}{\p \xi_i} = \alpha^{ij}_s,
\end{equation} 
where
\begin{align}\label{E:alphaijs}
    \alpha^{ij}_s = \left(\bar  Q^{ij}\right)^{(s-1)} - \bar\nabla^{dx^i} \left(u^j\right)^{(s-1)} + \bar\nabla^{dx^j} \left(u^i\right)^{(s-1)} + \\
 \frac{\p \pi^{ij}}{\p x^k} \left(u^k\right)^{(s-1)} + 
\sum_{t = 1}^{s-2}\{(u^i)^{(t+1)}, (u^j)^{(s-t)}\}_\Omega.\nonumber
\end{align}
Taking into account normalization condition (\ref{E:normcondloc}), we find a unique solution of (\ref{E:uniqueu}) which provides a recurrence formula for $u^i$:
\begin{equation}\label{E:reca}
    (u^i)^{(s)} = \frac{1}{s+1} \alpha^{ij}_s \xi_j.
\end{equation}

{\it Remark.} {\small The existence part of the proof of Theorem \ref{T:Fundamentalsolution} involves a simultaneous recursive calculation, along with $u^i$, of the local potentials $u_k$ and $v^i_k$ satisfying (\ref{E:vuagree}) and the condition
\[
    u^i = \pi^{ik}u_k.
\] 
This calculation is based upon the following observations. The function $\alpha^{ij}_s$ can be represented in the form
\begin{equation}\label{E:alphaaik}
   \alpha^{ij}_s = \pi^{ik} A^j_{s;k},
\end{equation}
where
\begin{align*}
   A^j_{s;k} =  (\bar A^j_k)^{(s-1)} - \frac{\p (u^j)^{(s-1)}}{\p x^k} + {^\dagger\!}\Gamma^{ml}_k \xi_l (v^j_m)^{(s-2)} + \pi^{jl} \frac{\p u_k^{(s-1)}}{\p x^l} + \\
 \Gamma^{jm}_l \xi_m \frac{\p u_k^{(s-1)}}{\p x^l} + \frac{\p \pi^{jl}}{\p x^k}u_l^{(s-1)} - \sum_{t = 1}^{s-2} \frac{\p (u^l)^{(t+1)}}{\p \xi_k} (v_l^j)^{(s-t+1)},
\end{align*}
and
\begin{align*}
\bar A^j_k =\frac{1}{2}\xi_p\xi_q\pi^{jm}\left( \frac{\p {^\dagger\!} \Gamma_m^{pq}}{\p x^k} 
 - \frac{\p {^\dagger\!} \Gamma_k^{pq}}{\p x^m}\right)  - 
\xi_p\xi_q {^\dagger\!}\Gamma^{lp}_k \Gamma^{jq}_l.
\end{align*}
This provides a recurrence formula for $u_k$:
\[
      u_k^{(s)} = \frac{1}{s+1} A^j_{s;k} \xi_j.
\]
On the one hand, to obtain a recurrence formula for $v_k^i$ we have to satisfy the condition
\[
   \frac{\p (u^i)^{(s)}}{\p \xi_j} = \pi^{jk} (v_k^i)^{(s-1)}.
\]
On the other hand, it follows from (\ref{E:reca}) that
\[
   \frac{\p (u^i)^{(s)}}{\p \xi_j} = \frac{1}{s+1}\left( \alpha^{ij}_s +  
\frac{\p \alpha^{il}_s}{\p \xi_j} \xi_l\right).
\]
Formulas (\ref{E:alphaaik}) and $\alpha^{ij}_s = -\alpha^{ji}_s$ imply that
\[
   \alpha^{ij}_s = \pi^{jk} (- A_{s;k}^i).
\]
One can check that
\[
   \frac{\p \alpha^{il}_s}{\p \xi_j}  =  \pi^{jk} B^{il}_{s;k},
\]
where
\begin{align*}
   B^{il}_{s;k} =  \left(\bar R^{il}_k\right)^{(s-2)} - \pi^{im} \frac{\p (v_k^l)^{(s-2)}}{\p x^m} + \pi^{lm} \frac{\p (v_k^i)^{(s-2)}}{\p x^m} - \\
\Gamma^{im}_n \xi_m \frac{\p (v_k^l)^{(s-2)}}{\p \xi_n} + \Gamma^{lm}_n \xi_m \frac{\p (v_k^i)^{(s-2)}}{\p \xi_n} +\\
 \Gamma^{im}_k (v_m^l)^{(s-2)} - \Gamma^{lm}_k (v_m^i)^{(s-2)} + \frac{\p \pi^{il}}{\p x^m}(v_k^m)^{(s-2)} +\\
\sum_{t=1}^{s-2} \left(- \frac{\p (v_k^i)^{(t)}}{\p \xi_m}(v_m^l)^{(s-t-1)} + (v_m^i)^{(t)}\frac{\p (v_k^l)^{(s-t-1)}}{\p \xi_m}\right)
\end{align*}
Thus we can give a recurrence formula for $v_k^i$:
\[
    (v_k^i)^{(s-1)} = \frac{1}{s+1}\left( - A_{s;k}^i +  B^{il}_{s;k}\xi_l\right).
\]

}
Next we want to give a recurrence formula for the unique global solution $F \in \H_\Omega$ of system (\ref{E:hamilt}) with the initial condition $F^{(0)} =f,\ f \in C^\infty(M)$.
Extracting from the equation
\[
     \bar\nabla^{dx^i} F = \{u^i,F\}_\Omega
\]
its homogeneous component of degree $s-1, s\geq 1,$ we obtain, using (\ref{E:firstterm}), the following formula:
\begin{equation}\label{E:fbeta}
    \frac{\p F^{(s)}}{\p \xi_i} = \beta^i_s,
\end{equation}
where
\begin{equation}\label{E:betais}
     \beta^i_s = \bar\nabla^{dx^i} F^{(s-1)} - \sum_{t = 1}^{s-1} \{(u^i)^{(t+1)}, F^{(s-t)}\}_\Omega.
\end{equation}
Thus we get a recurrence formula for $F$:
\begin{equation}\label{E:recb}
F^{(s)} = \frac{1}{s}\, \beta^i_s \xi_i.
\end{equation}

\section{Fedosov's Formal Symplectic Groupoids}

Assume that on a Poisson manifold $M$ with the Poisson tensor $\pi^{ij}$ there is a pair of associated global contravariant connections $(\nabla^\bullet, ^\dagger\! \nabla^\bullet)$ with $^\dagger\! \nabla^\bullet$ torsion-free, and on $(\T,Z)$ there is a global vertical 1-form $U$  with local potentials $u_k$ and $v^i_k$ as in Theorem \ref{T:lift}. These local potentials define associated local nonlinear invertible contravariant connections $D^\bullet$ and $^\dagger\! D^\bullet$ given by formulas (\ref{E:ncd}) and (\ref{E:ncdb}). Thus $^\dagger\! D^\bullet$ is torsion-free. Assume also that there exist global tensors $P^i_j$ and $Q^i_j$ on $M$ such that the contravariant connection $\nabla^\bullet$ respects them,
\begin{equation}\label{E:pq}
  P^i_j + Q^i_j = \delta^i_j, \mbox{ and } \pi^{ik} P^j_k = Q^i_k \pi^{kj}.
\end{equation}
We give two basic examples of the tensors $P^i_j$ and $Q^i_j$. The first one is $P^i_j = Q^i_j = (1/2) \delta^i_j$. The second example is related to the case of a K\"ahler-Poisson manifold endowed with the K\"ahler-Poisson contravariant connection (see Section \ref{S:linear}). In this case one can set
\begin{equation}\label{E:kahlpoisspq}
    P^k_m = \delta^k_m \mbox{ and } Q^{\bar l}_{\bar n} = \delta^{\bar l}_{\bar n},
\end{equation}
where $k,m$ are holomorphic and $\bar l, \bar n$ antiholomorphic indices. The tensor coefficients with mixed types of indices are set to zero. Otherwise speaking, $P$ projects a tangent vector onto its component of type (1,0) and $Q$ onto that of type (0,1).

Define a local change of variables on $(\T,Z),\ (x,\zeta) \mapsto (x, \xi)$ such that
\begin{equation}\label{E:xizeta}
   \xi_p = u_p\left(x, - \zeta_j P^j_\cdot\right) - u_p \left(x, \zeta_j Q^j_\cdot\right).
\end{equation}
We see from (\ref{E:invert}) and (\ref{E:pq}) that $\zeta_p = \xi_p \pmod{\xi^2}$, thus the change of variables is invertible. Introduce local mappings $f \mapsto Sf,\ f \mapsto Tf$ from functions on $M$ to functions on $(\T,Z)$ as follows:
\begin{equation}\label{E:sourcetarget}
   (Sf)(x,\xi) = \theta(f)\left(x, - \zeta_j P^j_\cdot\right)  \mbox{ and } 
    (Tf)(x,\xi) = \theta(f)\left(x, \zeta_j Q^j_\cdot\right).
\end{equation}
Denote by $\{\cdot,\cdot\}_\T$ the standard Poisson bracket on $\T$,
\[
    \{F,G\}_\T = \frac{\p F}{\p \xi_k}\frac{\p G}{\p x^k} - \frac{\p G}{\p \xi_k}\frac{\p F}{\p x^k}.
\]
It turns out that $S$ is a global Poisson morphism and $T$ is a global anti-Poisson morphism from $(C^\infty(M), \{\cdot,\cdot\})$ to $(C^\infty(\T,Z), \{\cdot,\cdot\}_\T)$, and for $f,g \in C^\infty(M)$ the elements $Sf$ and $Tg$ Poisson commute with respect to the Poisson bracket $\{\cdot,\cdot\}_\T$. Thus, according to \cite{FSG}, there exists a formal symplectic groupoid on $(\T,Z)$ whose source mapping is $S$ and target mapping is $T$. We call it {\it Fedosov's formal symplectic groupoid.}

One can prove the following lemma by straightforward calculations.
\begin{lemma}
  The space $\H_\Omega$ is closed with respect to the Poisson bracket $\{\cdot,\cdot\}_\T$.
\end{lemma}
 Thus $\H_\Omega$ has two different Poisson algebra structures corresponding to $\{\cdot,\cdot\}_\Omega$ and $\{\cdot,\cdot\}_\T$.
 We will prove that $S$ and $T$ are global mappings and that the images of $S$ and $T$ belong to $\H_\Omega$. The latter is a specific property of Fedosov's formal symplectic groupoids. First we need to prove a technical statement.
\begin{lemma}\label{L:matrixa}
  There exists an invertible formal matrix $A^i_k(x,\zeta)$ such that
\[
     \pi^{ik}\frac{\p \xi_k}{\p \zeta_l} = A^i_k \pi^{kl}.
\]
\end{lemma}
\begin{proof}
  Using formulas (\ref{E:vuagree}), (\ref{E:pq}), and (\ref{E:xizeta}), we proceed as follows:
\begin{align*}
   \pi^{ik}\frac{\p \xi_k}{\p \zeta_l} = -\pi^{ik} P^l_t \frac{\p u_k}{\p \xi_t}(x, -\zeta_j P^j_\cdot) - \pi^{ik} Q^l_t \frac{\p u_k}{\p \xi_t}(x, \zeta_j Q^j_\cdot) =\\
- \pi^{ts} P^l_t v^i_s (x, -\zeta_j P^j_\cdot) - \pi^{ts} Q^l_t v^i_s (x, \zeta_j Q^j_\cdot) = \\
\left( Q^s_k v^i_s (x, -\zeta_j P^j_\cdot) + P^s_k v^i_s (x, \zeta_j Q^j_\cdot)\right) \pi^{kl}.
\end{align*}
We see from (\ref{E:invert}) and (\ref{E:pq}) that the matrix
\[
   A^i_k = Q^s_k v^i_s (x, -\zeta_j P^j_\cdot) + P^s_k v^i_s (x, \zeta_j Q^j_\cdot)
\]
is invertible.
\end{proof}
Denote by $B^i_k$ the inverse of the matrix $A^i_k$. Lemma \ref{L:matrixa} implies the following formula:
\begin{equation}\label{E:bij}
    B^i_j \pi^{jk} = \pi^{il} \frac{\p \zeta_l}{\p \xi_k}.
\end{equation}
\begin{lemma}\label{L:sinhomega}
  For a given local function $f$ on $M$ there exists a formal function $a_l(x,\xi)$ such that
\[
   \frac{\p (Sf)}{\p \xi_k} = \pi^{kl} a_l.
\]
\end{lemma}
\begin{proof}
  It follows from Proposition \ref{P:poisson} that
\begin{equation}\label{E:dthetadxi}
   \frac{\p \theta(f)}{\p \xi_s} = \pi^{st}b_t
\end{equation}
for some formal function $b_t(x,\xi)$. Differentiating both sides of the first formula in (\ref{E:sourcetarget}) with respect to $\zeta_l$ we get, using (\ref{E:dthetadxi}), the following equality:
\begin{equation}\label{E:equality}
   \frac{\p (Sf)}{\p \xi_k} \frac{\p \xi_k}{\p \zeta_l} = - P^l_s \pi^{st}b_t \left(x, - \zeta_j P^j_\cdot\right). 
\end{equation}
Taking into account (\ref{E:pq}), (\ref{E:bij}), and (\ref{E:dthetadxi}), we get from (\ref{E:equality}) that
\[
    \frac{\p (Sf)}{\p \xi_k} = - \pi^{kl} B^s_l Q^t_s b_t \left(x, - \zeta_j P^j_\cdot\right) 
\]
whence the lemma follows.
\end{proof}
\begin{lemma}\label{L:auxil}
Assume that $F(x,\xi)$ is a formal function such that
\[
   \frac{\p F}{\p \xi_k} = \pi^{kl} a_l
\]
and $\phi_p(x,\xi),\psi_p(x,\xi)$ are formal functions such that $\pi^{ip}\phi_p = \pi^{ip}\psi_p$. Then $F(x,\phi) = F(x,\psi)$.
\end{lemma}
\begin{proof}
Consider the function
\[
    \alpha(t) = F(x, t \phi + (1-t) \psi),
\]
so that $\alpha(0) = F(x,\psi)$ and $\alpha(1) = F(x,\phi)$. We have
\[
    \frac{d\alpha}{dt} = \frac{\p F}{\p \xi_k} (\phi_k - \psi_k) = \pi^{kl} a_l (\phi_k - \psi_k) =0,
\]
whence the lemma follows.
\end{proof}
The local potential $u_p$ of the global vertical 1-form $U$ is not uniquely defined. Assume that $\tilde u_p$ is another local potential of $U$. Then $\pi^{ip} u_p = \pi^{ip}\tilde u_p$. It follows from lemmas \ref{L:sinhomega} and \ref{L:auxil} that 
formulas (\ref{E:xizeta}) and (\ref{E:sourcetarget}) determine the same mapping $S$ for the potentials $u_p$ and $\tilde u_p$. Thus $S$ is globally defined. Lemma \ref{L:sinhomega} implies that the image of $S$ is in the space $\H_\Omega$.
The proof that $T$ is globally defined and has the image in $\H_\Omega$ is similar.
The following theorem can be proved by long straightforward calculations.
\begin{theorem}
  The mappings 
\[
S,T : (C^\infty(M),\{\cdot,\cdot\}) \to (\H_\Omega, \{\cdot,\cdot\}_\T)
\]
are a Poisson an anti-Poisson morphisms, respectively. Moreover, the images of the mappings $S$ and $T$ Poisson commute. Thus there is a unique formal symplectic groupoid on $(\T,Z)$ whose source and target maps are $S$ and $T$, respectively.
\end{theorem}

\section{The symplectic case}

In \cite{Deq} we gave a self-contained construction of the formal symplectic groupoid
of Fedosov's deformation quantization. In this section we will show that in the symplectic case Fedosov's formal symplectic groupoid constructed in the previous sections is the same groupoid as in \cite{Deq}. Recall the construction from \cite{Deq}. 

Let $M$ be a symplectic manifold with symplectic form $\omega = \frac{1}{2}\omega_{ij} dx^i \wedge dx^j$. The inverse of $\omega_{ij}$ is a Poisson tensor $\pi^{ij}$. We assume that there is a tensor $\Lambda^{ij}$ on $M$ such that
\[
    \frac{1}{2}\left(\Lambda^{ij} - \Lambda^{ji}\right) = \pi^{ij}.
\]
Let $\nabla_\bullet$ be a covariant connection which respects $\Lambda^{ij}$ (and therefore the symplectic and Poisson tensors as well). The connection $\nabla_\bullet$ may have torsion. From these data one can construct Fedosov's deformation quantization on $M$ and the corresponding formal symplectic groupoid (see \cite{Deq},\cite{FSG}). The groupoid construction starts with the ($\nu$-free) Fedosov's lift 
\[
\tau^\vee: C^\infty(M) \to C^\infty(TM,Z)
\]
first extracted from Fedosov's quantization in \cite{EW}.

Introduce Fedosov's operators on $C^\infty(TM)\otimes \Omega^1(M)$:
\[
    \delta (a) = dx^i \wedge \frac{\p a}{\p y^i}, \quad \delta'(a) = 
    y^i\, \iota\left(\frac{\p}{\p x^i}\right) a,
\]
where $\{x^i\}$ are local coordinates on $M$ and $\{y^i\}$ are the corresponding fibre coordinates on $TM$. Denote by $\{\cdot,\cdot\}_{TM}$ the fibrewise Poisson bracket on $TM$,
\[
    \{F,G\}_{TM} = \pi^{jk} \frac{\p F}{\p y^j} \frac{\p G}{\p y^k}.
\]
Connection $\nabla_\bullet$ on $M$ induces a covariant connection $\bar \nabla_\bullet$ on $(TM,Z)$ expressed locally as follows:
\[
   \bar\nabla_{dx^i} = \frac{\p}{\p x^i} - \Gamma_{ij}^k y^j \frac{\p}{\p y^k}.
\]
 There exists a flat nonlinear covariant connection $D_\bullet$ on $(TM,Z)$,
\[
   D_\bullet = \bar\nabla_\bullet + \{\rho,\cdot\}_{TM},
\]
where $\rho = \rho_p(x,y) dx^p$ is a unique 1-form satisfying the equation
\[
   - \omega = R + \bar\nabla_\bullet \rho + \frac{1}{2}\{\rho,\rho\}_{TM}
\]
and such that $\delta' \rho = 0$ and $\rho_p = y^k \omega_{kp} \pmod{y^2}$. Here
\[
  R = \frac{1}{4}\omega_{s\alpha} R^\alpha_{tkl} y^s y^t dx^k \wedge dx^l,
\]
where 
\[
    R^s_{tkl} = \frac{\p \Gamma^s_{lt}}{\p x^k} - \frac{\p \Gamma^s_{kt}}{\p x^l} + \Gamma^s_{k\alpha}\Gamma^\alpha_{lt} - \Gamma^s_{l\alpha}\Gamma^\alpha_{kt}
\]
is the curvature of $\nabla_\bullet$. In the notations of \cite{Deq},
\[
   \rho_p = y^k \omega_{kp} + r_p^\vee.
\]
Fedosov's lift $F = \tau^\vee(f)$ of a function $f \in C^\infty(M)$ is a unique element of $C^\infty(TM,Z)$ such that $D_\bullet F = 0$ and $F|_{y = 0} = f$. The mapping $\tau^\vee$ is a Poisson morphism from $C^\infty(M)$ to $(C^\infty(TM,Z), \{\cdot,\cdot\}_{TM})$.
Consider the global diffeomorphism of $(\T,Z), \ (x,\zeta) \mapsto (x, \xi),$ such that
\[
   \xi_p = \rho_p\left(x,\frac{1}{2}\Lambda^{j\cdot} \zeta_j\right) - \rho_p\left(x, \frac{1}{2}\Lambda^{\cdot j} \zeta_j\right).
\]
Then, according to \cite{Deq}, the source and the target mappings of the formal symplectic groupoid of Fedosov's deformation quantization are given by the formulas
\[
   (Sf)(x,\xi) = \tau^\vee (f) \left(x, \frac{1}{2}\Lambda^{\cdot j} \zeta_j\right) \mbox{ and } (Tf)(x,\xi) = \tau^\vee(f) \left(x, \frac{1}{2}\Lambda^{j\cdot} \zeta_j\right),
\] 
where $f \in C^\infty(M)$. To obtain the construction of Fedosov's formal symplectic groupoid one has to pull back Fedosov's lift via the mapping $\#$. For a function $F = F(x,y) \in C^\infty(TM,Z)$ denote its pullback via $\#$ by $F^\#$ so that 
\[
          F^\#(x, \xi) = F(x, \xi_j \pi^{j\cdot}).
\]
We get the following identifications. For $f \in C^\infty(M)$
\[
   \theta(f) = \left(\tau^\vee (f)\right)^\#, u_p = - \rho_p^\#, \mbox{ and } v^q_p = - \left(\pi^{qs} \frac{\p \rho_s}{\p y^p}\right)^\#.
\]
Finally,
\[
    P^j_k = \frac{1}{2} \omega_{ki} \Lambda^{ij} \mbox{ and } Q^j_k = \frac{1}{2} \Lambda^{ji} \omega_{ik}.
\]
Notice that in the symplectic case $u_p dx^p$ is a global object on $M$, while in the general Poisson case only $\pi^{pq} u_q \frac{\p}{\p x^p}$ is global.

\section{The K\"ahler-Poisson case}

Assume that $M$ is a K\"ahler-Poisson manifold with the K\"ahler-Poisson tensor 
$g^{\bar l k}$, where $k, \bar l$ are holomorphic and antiholomorphic indices, respectively (we use the agreement that $\pi^{\bar lk}= g^{\bar lk},\ \pi^{k\bar l}= - g^{\bar lk}, \ \pi^{km} =0$ and $\pi^{\bar l\bar n} =0$).
Denote by $\nabla^\bullet$ the K\"ahler-Poisson contravariant connection whose Christoffel symbols are given by (\ref{E:christoff}). The K\"ahler-Poisson contravariant connection is Poisson and torsion-free. Introduce tensors $P$ and $Q$ on $M$ by formula (\ref{E:kahlpoisspq}). With these data, one can construct a Fedosov's formal symplectic groupoid on $(\T,Z)$. It was proved in \cite{FSG} that for a K\"ahler-Poisson manifold $M$ there exists a unique formal symplectic groupoid with separation of variables on $(\T,Z)$, i.e., such that for a local holomorphic function $a$ and a local antiholomorphic function $b$ on $M$
\begin{equation}\label{E:fsgsv}
    Sa =a \mbox{ and } Tb = b.
\end{equation}
We want to prove the following proposition.

\begin{proposition}\label{P:kahlpoiss}
  The Fedosov's formal symplectic groupoid on $(\T,Z)$ constructed from the data $(\nabla^\bullet, P,Q)$ is the formal symplectic groupoid with separation of variables.
\end{proposition}

For local holomorphic coordinates $z^k,\bar z^l$ on $M$ denote by $\eta_k,\bar \eta_l$ the corresponding fibre coordinates on $\T$. Denote by $\I$ and $\bar \I$ the ideals in $C^\infty(\T,Z)$ locally generated by the coordinates $\eta$ and $\bar \eta$, respectively.
Formulas
\[
   \bar\nabla^{dz^k} = - g^{\bar lk} \frac{\p}{\p \bar z^l} + \frac{\p g^{\bar n k}}{\p \bar z^ l} \bar \eta_n \frac{\p}{\p \bar \eta_l}, \ \bar\nabla^{d\bar z^l} = g^{\bar lk} \frac{\p}{\p z^k} - \frac{\p g^{\bar l m}}{\p z^k} \eta_m \frac{\p}{\p \eta_k} 
\]
imply that the connection $\bar\nabla^\bullet$ respects the ideals $\I$ and $\bar \I$.
Let us show that the spaces $\H_\Omega \cap \I$ and $\H_\Omega \cap \bar\I$ are closed with respect to the Poisson bracket $\{\cdot,\cdot\}_\Omega$. Assume that $F,G \in \H_\Omega \cap \I$. Then there exist local potentials $a_k, a_{\bar l}$ and $b_k,b_{\bar l}$ of $F$ and $G$, respectively, such that
\[
   \frac{\p F}{\p \bar \eta_l} = g^{\bar lk} a_k, \ \frac{\p F}{\p \eta_k} = - g^{\bar lk} a_{\bar l}, \ \frac{\p G}{\p \bar \eta_l} = g^{\bar lk} b_k, \ \frac{\p G}{\p \eta_k} = - g^{\bar lk} b_{\bar l}.
\]
Now,
\[
   \{F,G\}_\Omega = g^{\bar lk}a_{\bar l} b_k  - g^{\bar lk}a_k b_{\bar l} = \frac{\p G}{\p \bar \eta_l}a_{\bar l} - \frac{\p F}{\p \bar \eta_l} b_{\bar l} \in \I.
\]
The proof that $\H_\Omega \cap\bar\I$ is closed with respect to the Poisson bracket $\{\cdot,\cdot\}_\Omega$ is similar. Using these statements and recurrence formulas (\ref{E:reca}) and (\ref{E:recb}) one can prove by induction the following lemma.
\begin{lemma}\label{L:ufetabareta}
 (a) Let $U$ be the normalized solution of the Fundamental Equation corresponding to the 
K\"ahler-Poisson connection $\bar\nabla^\bullet$ on $M$. Then its local potential $u_k, u_{\bar l}$ satisfies the condition
\[
   u_k = - \eta_k, u_{\bar l} = - \eta_{\bar l} \pmod{\I \cap \bar \I}.
\]
(b) Let $a$ and $b$ be a local holomorphic and a local antiholomorphic functions on $M$, respectively. Then
\[
   \theta(a) - a \in \bar\I \mbox{ and } \theta(b) - b \in \I. 
\]
\end{lemma}
Lemma \ref{L:ufetabareta} implies that for the tensors $P$ and $Q$ given by formula (\ref{E:kahlpoisspq}) change of variables (\ref{E:xizeta}) is trivial. Namely, for $u_k = u_k(z,\bar z, \eta, \bar\eta)$ and $u_{\bar l} = u_{\bar l}(z,\bar z, \eta, \bar\eta)$
\begin{align}\label{E:trivchange}
    u_k(z,\bar z, - \eta,0) - u_k(z,\bar z, 0, \bar \eta) = \eta_k \mbox{ and }\\
  u_{\bar l}(z,\bar z, - \eta,0) - u_{\bar l}(z,\bar z, 0, \bar \eta) = \eta_{\bar l}.\nonumber 
\end{align}
Now conditions (\ref{E:fsgsv}) follow from formulas (\ref{E:sourcetarget}), (\ref{E:trivchange}), and Lemma \ref{L:ufetabareta}, which proves Proposition \ref{P:kahlpoiss}.

\end{document}